\documentclass[sn-mathphys-num]{sn-jnl}

\usepackage{graphicx}%
\usepackage{multirow}%
\usepackage{amsmath,amssymb,amsfonts}%
\usepackage{amsthm}%
\usepackage{mathrsfs}%
\usepackage[title]{appendix}%
\usepackage{xcolor}%
\usepackage{textcomp}%
\usepackage{manyfoot}%
\usepackage{booktabs}%
\usepackage{threeparttable}
\usepackage{algorithm}%
\usepackage{algorithmicx}%
\usepackage{algpseudocode}%
\usepackage{listings}%
\usepackage{subcaption}

\newcommand*{\img}[1]{%
    \raisebox{-.3\baselineskip}{%
        \includegraphics[
        height=\baselineskip,
        width=\baselineskip,
        keepaspectratio,
        ]{#1}%
    }%
}

\theoremstyle{thmstyleone}%
\theoremstyle{thmstyletwo}%

\theoremstyle{thmstylethree}%

\begin{document}

\title[Multi-timescale Stochastic Programming with Applications in Power Systems]{Multi-timescale Stochastic Programming with Applications in Power Systems}

\author[1]{\fnm{Yihang} \sur{Zhang}}\email{zhangyih@usc.edu}

\author[1]{\fnm{Suvrajeet} \sur{Sen}}\email{suvrajes@usc.edu}

\affil[1]{\orgdiv{Daniel J. Epstein Department of Industrial \& Systems Engineering}, \orgname{University of Southern California}, \orgaddress{\street{3650 McClintock Ave}, \city{Los Angeles}, \postcode{90089}, \state{CA}, \country{USA}}}


\abstract{This paper introduces a multi-timescale stochastic programming framework designed to address decision-making challenges in power systems, particularly those with high renewable energy penetration. The framework models interactions across different timescales using aggregated state variables to coordinate decisions. In addition to Multi-timescale uncertainty modeled via multihorizon trees, we also introduce a "synchronized state approximation," which periodically aligns states across timescales to maintain consistency and tractability. Using this approximation, we propose two instantiation methods: a scenario-based approach and a value function-based approach specialized for this setup. Our framework is very generic, and covers a wide-spectrum of applications. }




\maketitle

\section{Introduction}

Renewable energy (e.g., wind, solar, etc.) is expected to supply a significant portion of power to the electricity grid. These resources help lower the carbon footprint of an energy system because they can replace energy produced by fossil fuel-based power plants. However, the variable nature of renewable resources poses its own challenges for electric system operators, who are required to meet reliability standards set by utility commissions. In addition, diurnal variations in the atmosphere are known to cause further mismatches between supply and demand (for electricity). Such mismatches (e.g., the so-called ``duck-chart'' phenomenon observed in California) can, in fact, amplify the challenges associated with ensuring system reliability. 
Due to the influx of renewable energy, the Federal Regulatory Commission (FERC) has mandated \cite{federal_energy_regulatory_commission_integration_2012} that time discretization of electricity dispatch planning must be reduced to maintain system reliability. In fact, the recommendation is to reduce the planning interval from the current 60-minutes to much shorter (e.g., 15-minute) intervals. Even with such a reduction in time steps, one must accommodate wind variability, whose speed and direction can change in very short time-intervals (e.g., a couple of minutes). The implication of these changes is that traditional deterministic approaches to power system planning is difficult to justify. The system, ideally, needs to have at least a long enough decision horizon to recognize trends, as well as a short enough time-resolution to represent environmental impact well enough, so that the carbon-footprint can be reduced, without sacrificing system reliability.  The challenges are clearly multi-faceted due to the interactions between technology, market economics, and operational constraints, including computational challenges of developing reasonable control strategies at scale.  Because long run sustainability of the planet requires a multi-faceted outlook, it is important for energy planning models to accommodate multiple factors at different scales.  These factors have already been acknowledged in recent industry-driven publications, although the associated science of decision-making (models and algorithms) have yet to be designed and tested. 

As the reader will observe from this presentation, Multi-timescale stochastic programming models must overcome many challenges because of several considerations (e.g., dependence across scales), and the interplay between data and decisions in a manner which respects the multi-scale time-dependent evolution of data and nonanticipativity of decisions. For example, using traditional stochastic programming methodology, modeling a 24-hour horizon at a resolution of 15 min would, perhaps, require 96 stages.  Constructing a reasonable scenario tree for such a large number of stages is nearly impossible, let alone solving the resulting problem using standard off-the-shelf solvers. We should emphasize that the idea of non-anticipativity itself is a complicated phenomenon in multi-timescale models, and highlights the need for specialized approaches tailored to multi-timescale models.

\textbf{Our Contribution} In this paper we outline a mathematically precise statement of a general-purpose multi-timescale stochastic program. While approximations have been set up for some multi-timescale problems (e.g., using multi-horizon scenario trees as in \cite{kaut_multi-horizon_2014}), they are intended to leverage standard stochastic programming methodologies. Our goal is to create a scalable model which has the potential to be resolved using modern tools of stochastic programming, e.g., adaptive sampling, that is reminiscent of a combination of stochastic programming and dynamic programming.

The scales in our problem refer to \textit{timescales}. This is conceptually different from methods that deal with joint modeling and optimization at various \textit{spatial} scales, as seen in chemical engineering, and material science, for instance, \cite{biegler_multi-scale_2014,miller_next_2018}. Although multi-spatial-scale methods exist in power systems, such as in microgrid applications \cite{carpentier_upper_2019}, and while they can (and perhaps should) be used together with multi-timescale methods, such integration is beyond the scope of this paper.

\textbf{Outline} The current paper is structured as follows: We begin in Section \ref{sec:motivational_example} with a motivational example that highlights the key challenges of multi-timescale decision-making. In Section \ref{section:comparison_with_existing_methods}, we then review existing modeling and algorithmic approaches from the literature.  In Section \ref{section:formulation}, we define the time notation and the key elements of a multi-timescale stochastic program. Importantly, we define the concept of "aggregated state variables," which serve as the communicative link between different timescales. Section \ref{section:equivalent_multistage_problem} introduces the formulation of a multi-timescale problem as an equivalent multistage problem. The resulting problem is unlikely to be tractable for off-the-shelf solvers. As a remedy, we propose in Section \ref{section:tractability} an approximation called \textbf{synchronized state approximation} to use with \cite{kaut_multi-horizon_2014} and \cite{maggioni_bounds_2020}. Then, in Section \ref{section:instantiation}, we show how to instantiate the problem from the description: either we do a scenario-based instantiation or realize these using a type of special value function. Finally, in Section \ref{sec:multi-timescale_instance}, we present a detailed three-timescale model of the motivational problem. This example demonstrates a hybrid approach, mixing scenario-based and value function-based instantiations to best leverage the characteristics of the decisions and data processes at each timescale.

\section{A Motivational Example} \label{sec:motivational_example}
A motivational example provided by Atakan et al. \cite{atakan_towards_2022} introduced three timescales. They are called \textbf{day ahead (DA)}, represented by blue squares in figures (\img{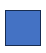}), \textbf{short term (ST)}, represented by green triangles (\img{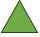}), and \textbf{real time (RT)}, represented by yellow circles (\img{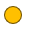}).

\begin{table}[ht]
\centering
\begin{tabular}{|p{2cm}|p{2.5cm}|p{3cm}|p{3cm}|}
\hline
 & \textbf{Day Ahead} \img{square.PNG} & \textbf{Short Term} \img{triangle.PNG} & \textbf{Real Time} \img{circle.PNG} \\ \hline
\textbf{Decision time} & 0:00, 4:00, 8:00, \dots & 0:00, 1:00, 2:00, \dots & 0:00, 0:15, 0:30, 0:45, 1:00, \dots \\ \hline
\textbf{Systems} & Slow-ramping generators & Fast-ramping generators & Renewable generators and batteries \\ \hline
\textbf{States} & \multicolumn{2}{|c|}{Operation status, power output}  & Batteries charge  \\ \hline
\textbf{Control} & \multicolumn{2}{|c|}{Turn on/off, ramp up/down} & Power flow, charge, discharge batteries  \\ \hline
\textbf{Stochastics} & Autocorrelated prediction error & Markovian switching & White noise residuals \\ \hline
\end{tabular}
\caption{Multiple timescale modeling of a power system from \cite{atakan_towards_2022}}
\label{tab:generators_dispatch}
\end{table}

Decisions within the Day Ahead timescale are made every four hours, specifically at 0:00, 4:00, and 8:00. During these intervals, operators can ramp up, ramp down, or turn slow-ramping generators on or off. Short-term decisions are made hourly, at 0:00, 1:00, and 2:00 ... These decisions manage the fast-ramping generators. The Real-Time timescale operates at a higher granular level of granularity, with decisions made every 15 minutes, at times such as 0:00, 0:15, 0:30, etc. At this scale, one can observe the actual renewable generation. Battery storage may be charged or discharged as needed. Load shedding and curtailment may also be realized, with ``penalty" costs applied accordingly. Note that while the power output from renewable generators is not directly controllable, the operator can decide whether to use the generated power.

The stochastics capture the error between the prediction and actual realizations of the renewable production, inspired by the model proposed in \cite{goujard_mape_maker_2023}. In the Day-Ahead timescale, the autocorrelated prediction error is modeled, and it is possible to model stochastic evolution using scenarios. In the short-term timescale (hourly), we can use Markov-switching to represent finer adjustment to the prediction. Finally, the Real-Time timescale describes the residual prediction error (see \cite{cao_switching_2016}) that is not captured by coarser timescales. However, here, we simply assume these residual errors are white-noise.

\section{Review of existing methods} \label{section:comparison_with_existing_methods}
\subsection{Direct Stochastic Programming Approaches}
Below are a few multistage algorithms that can be used to solve the equivalent multistage problem of the multi-timescale problem directly.
\begin{itemize}
    \item Price-directive decomposition methods (progressive hedging and its variants, e.g., randomized\cite{bareilles_randomized_2020}, integer heuristics\cite{watson_progressive_2011}, integer\cite{atakan_progressive_2018}) are effective at solving stochastic unit commitment problems. The solver can take advantage of the tight formulation of the generator polyhedron in the subproblems. However, as stated in the previous sections, such direct application does not scale up well because of the exponential explosion in the number of scenarios.
    \item Resource-directive decomposition methods (e.g., stochastic dual dynamic programming SDDP and its integer variant SDDiP\cite{ahmed_stochastic_2016}) provide a compact description of the problem using dynamic programming, which greatly increases tractability. The drawback is that SDDiP requires pure binary state variables, while typically, generator polytopes require mixed-binary state variables, so SDDiP is not guaranteed to converge when applied to the unit commitment problem without binary expansion in the continuous variables. It is also difficult to model long-term non-markovian effects such as minimum uptime/downtime and correlations.
    \item Dynamic programming methods are also good for small-scale problems and naturally allow general integers. They are not scalable because of the combinatorial nature of the state variables.
\end{itemize}

\subsection{Multi-horizon Tree Methods \cite{kaut_multi-horizon_2014}}

Multi-horizon tree is one of the first attempts to model this multi-timescale problem. Recently, algorithms like \cite{zhang_stabilised_2024} have been developed to solve these multi-horizon stochastic programs. We believe there are two major differences between our modeling approach and the multi-horizon approach.

First, our work is focused on the coordination between timescales and maintaining a feasible trajectory in all timescales across the whole horizon, whereas, in a multi-horizon tree, it is only required to maintain a valid state trajectory in each ``segment'' problem. In some applications, such as the unit commitment problem, the requirement for a coherent state trajectory cannot be reasonably discarded.

Secondly, the multi-horizon tree paper focuses on the scenario-based formulation, while our formulation can be modified to use some resource-directive decomposition methods.

\subsection{Other Multi-timescale Modeling Approaches}
\begin{figure}
    \centering
    \includegraphics[width=\linewidth]{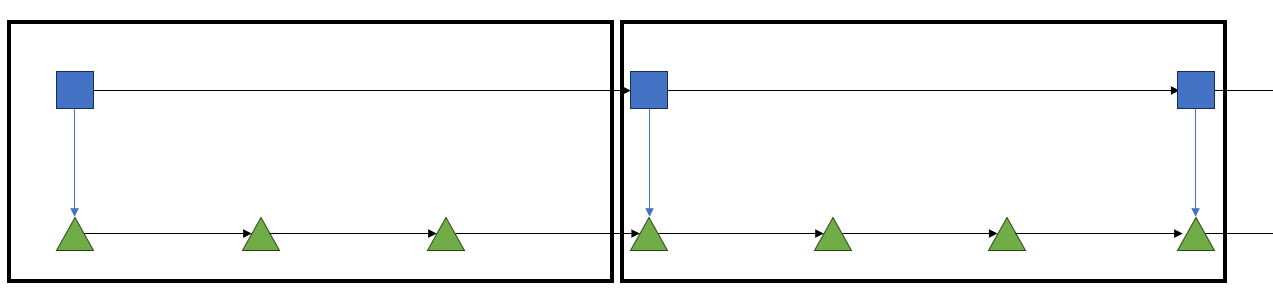}
    \caption{Illustration of Blockwise Independence Structure in \cite{carpentier_time_2023}, \cite{rigaut_decomposition_2024}}
        \label{fig:blockwise_independence}
\end{figure}

\textbf{Carpenter and De Lara \cite{carpentier_time_2023}, \cite{rigaut_decomposition_2024} }
The line of work by Carpenter and De Lara describes a similar situation where the decisions are made on multiple scales. In their work, they used a 'white noise' assumption on the slow timescale, which is a blockwise independence assumption. With some effort, this could be relaxed to a blockwise Markovian assumption. Such a simplifying assumption can be justified in many applications. By exploiting this independence structure, they are able to form a temporal decomposition scheme along the slow timescale, which greatly increases the tractability of the problem. Our method is more focused on scenario trees, so it is less tractable, but we are able to model dependence structure in the randomness, specifically, the autocorrelatedness of the renewable generation in the slow timescale. Figure \ref{fig:blockwise_independence} shows the independence assumption used in their work.

\textbf{Glanzer and Pflug\cite{glanzer_multiscale_2020}}
Pflug et al. describe a multi-timescale problem. In their work, they used scenario lattice and stochastic bridges to describe the stochastics in a consistent way across the timescales. However, their work mainly focused on the modeling side but did not address a particular algorithm to solve the resulting optimization problem. The modeling technique in their paper can be applied and realized in an efficient way in our framework using a specialized value function, described in the later section.

\textbf{Gangammanavar et al. \cite{gangammanavar_stochastic_2016}, Atakan et al. \cite{atakan_towards_2022}} These works state the unit commitment problem in multiple timescales. They apply a rolling horizon scheme by sequentially solving two-stage stochastic programs as an approximation and rolling the horizon forward by fixing some variables to the optimal solution of the two-stage approximations. This method can be viewed as a sequential approximation scheme of our multi-timescale problem.

\section{Formulation of a Multi-timescale Problem} \label{section:formulation}

In this section, we introduce the notation used in a multi-timescale stochastic program. It is necessary to use functional notation for both state and control variables, implying that these variables are understood to be random variables themselves. In order to help readers obtain a clearer sense what this notation implies, the Appendix includes section \ref{section:appendix-multistage} where we present a single-timescale multistage problem using this notation.

\subsection{Timescales and Ticks} \label{subsection:multitimescale_formulation}

We represent the decisions point on different timescales with a set of ticks and assign a zero-based tuple to index each tick. At each tick, the decision maker is required to produce a decision on that corresponding timescale. The ticks are analogous to the notion of \textit{stages} in a single-timescale model.

The arrangement of the ticks in the motivational example is summarized in Table \ref{table:ticks}.

\begin{table}[] 
\centering
\begin{tabular}{@{}lllllllllll@{}}
\toprule
        & 0:00  & 0:15  & … & 1:00  & 1:15  & … & 2:00  & … & 4:00  & 4:15  \\ \midrule
DA \includegraphics[scale=0.3]{square.PNG} & 0     &       &   &       &       &   &       &   & 1     &       \\
ST \includegraphics[scale=0.3]{triangle.PNG} & 0,0   &       &   & 0,1   &       &   & 0,2   & … & 1,0   &       \\
RT \includegraphics[scale=0.5]{circle.PNG} & 0,0,0 & 0,0,1 & … & 0,1,0 & 0,1,1 & … & 0,2,0 & … & 1,0,0 & 1,0,1 \\ \bottomrule
\end{tabular}
\caption{Indexing Timescales in the motivational example}
\label{table:ticks}
\end{table}

The DA ticks (in the slowest timescale) are indexed by 1-tuple, e.g., $(0), (1), (2), \dots $, which refers to the tick at 0:00 a.m., 4:00 a.m., 8:00 a.m, etc., in our motivational example. The ticks in the next timescale layer are indexed by 2-tuples $(0,0), (0,1), (0,2), (0,3), (1,0), (1,1), \dots$, meaning 0:00, 1:00, 2:00, 3:00, 4:00, 5:00 and so on. Collectively, we define the set of ticks in each timescale as $\mathcal T_1 = \{(0), (1), (2), \dots\}, \mathcal T_2 = \{(0,0), (0,1), \dots\}, \mathcal T_3 = \{(0,0,0), (0,0,1), \dots \}$. To use this notation, we need to ensure that the ticks on the slower timescales coincide with those on faster timescales.

For a tick $\tau_j \in \mathcal T_j$, we use $\tau_j+$ to mean \textbf{the next tick on the same timescale}. For example, the next tick of the DA tick 0 (at 0:00) is tick 1 (at 4:00); the next tick of ST tick (0,1) (at 1:00) is tick (0,2) (at 2:00). The \textbf{parent tick} $\mathcal T_{j-1}(\tau_{j})$ means the tick in the $j-1$-th timescale that contains the $j$-th timescale tick. For example, the parent tick of (0,1,1) is (0,1) and so on.

In the example above, although the ticks on each timescale are set up at uniform intervals, it is not a requirement of the multiscale model formulation. We can adopt a dynamic resolution for periods requiring more precise modeling, for example, during the duck neck part of the demand curve.

\subsection{Multi-timescale model without interaction between timescales}

Below, we show the components of a multi-timescale problem along with its meaning for our motivational problem. This is a three-timescale example, and the same definition can be extended to more timescales.

\begin{itemize}
    \item Ticks $\mathcal T_1 = \{0, 1, 2 \dots\}, \mathcal T_2 = \{(0,0), (0,1), \dots\}, \mathcal T_3 = \{ (0,0,0), \dots \}$.
\end{itemize}

For each tick \(\tau_j \in \mathcal{T}_j\), where \(j = 1, 2, 3\):

\begin{itemize}
    \item State $x_{\tau_j}$ (meaning whether a generator is on/off, current power output at $\tau_j$.) \footnote{This may also include information about uptime or downtime of a generator if necessary.}
    \item Control $u_{\tau_j} \in U_{\tau_j}(x_{\tau_j})$ (whether we would like to turn on/off, ramp up/down the generator)
    \item Randomness $\omega_{\tau_j} \in \Omega_{\tau_j}$ (timescale dependent randomness, the renewable generation at $\tau_j$)
    \item State Dynamics $x_{{\tau_j}+} = D_{\tau_j}(x_{\tau_j}, u_{\tau_j}, \omega_{\tau_j})$
    \item (Instantaneous) Cost $f_{\tau_j}(x_{\tau_j}, u_{\tau_j})$ 
\end{itemize}

The within-timescale dynamics are defined in the same manner as in a single-timescale problem. It is worth emphasizing again that both \(x_{\tau_j}\) and \(u_{\tau_j}\) are random decision variables as functionals.

Putting the above together, we have the following multi-timescale problem.

\begin{align}
    \text{Min } & E[\sum_j \sum_{\tau_j \in \mathcal T_j} f(x_{\tau_j}, u_{\tau_j})] \\
    & x_{{\tau_j}+} = D_{\tau_j}(x_{\tau_j}, u_{\tau_j}, \omega_{\tau_j}) \\ \label{eqn:multitimescale_no_adj_state}
    & u_{\tau_j} \in U_{\tau_j}(x_{\tau_j}) \qquad \forall \tau_j \in \mathcal T_j\\ \label{eqn:multitimescale_no_adj_action}
    & (x_{\tau_j}, u_{\tau_j}) \text{ non-anticipative} 
\end{align}

In some applications, we need to model history-dependent state dynamics, where the state transition depends on all previous states. This can be done by replacing \eqref{eqn:multitimescale_no_adj_state} with the following equation.

\begin{equation}
    x_{{\tau_j}+} = D_{\tau_j}(\{x_{\tau}\}_{\tau \in \mathcal T_j: \tau \leq \tau_j}, u_{\tau_j}, \omega_{\tau_j})
\end{equation}

We defer the detailed discussion of non-anticipativity to Section \ref{section:nac_multitimescale}.

\subsection{Modeling Interaction Between Timescales Using Aggregated State Variables}

We observe that in the previously defined multi-timescale problem, there is no interaction whatsoever between the systems operating at different timescales, so the problem naturally decomposes into different timescales as individual multistage problems.

The interaction between the decisions in different timescales is why this model is useful. In the motivational example, a tick in the Day-Ahead (DA) timescale aggregates all available generator power and reserves committed during that tick. This information needs to be conveyed to the ticks in the Short-Term (ST) timescale so these faster ST ticks can make informed decisions based on the summarized information from the parent tick. Conversely, if necessary, the ticks in the ST timescale can give feedback to the DA layer, allowing adjustments based on developments in the faster timescales. This dynamic interaction also applies between the ST and the Real Time (RT) timescales in the same way.

To model this interaction, we introduce \textbf{aggregated state variables}, denoted as \(y_{\tau}\), for each tick \(\tau_j\) except for the fastest timescale. These variables function similarly to state variables but are specifically designed to coordinate decisions between timescales.

In each tick (except those in the fastest timescale), the aggregated state variable is generated from the state variable by adding the following constraint.
\begin{equation} \label{eqn:agg_state_generation}
    y_{{\tau_j}} = C_{\tau_j}(x_{\tau_j}, \omega_{\tau_j})
\end{equation}

Correspondingly, in each fast timescale $\tau_j$, we let feasible controls and state dynamics depend on the aggregated state they receive from their parent ticks $\mathcal{T}_{j-1}(\tau_j)$. \eqref{eqn:multitimescale_no_adj_state} and \eqref{eqn:multitimescale_no_adj_action} should be adjusted as follows:

\begin{itemize}
    \item Control 
    \begin{equation} \label{eqn:multiscale_feasible_actions}
        u_{\tau_j} \in U_{\tau_j}(x_{\tau_j}, y_{\mathcal{T}_{j-1}(\tau_j)})
    \end{equation}
    \item State Dynamics
    \begin{equation} \label{eqn:multiscale_transition}
        x_{{\tau_j}+} = D_{\tau_j}(x_{\tau_j}, u_{\tau_j}, y_{\mathcal{T}_{j-1}(\tau_j)}, \omega_{\tau_j})
    \end{equation}
\end{itemize}

\begin{figure}
    \centering
     \includegraphics[width=\linewidth]{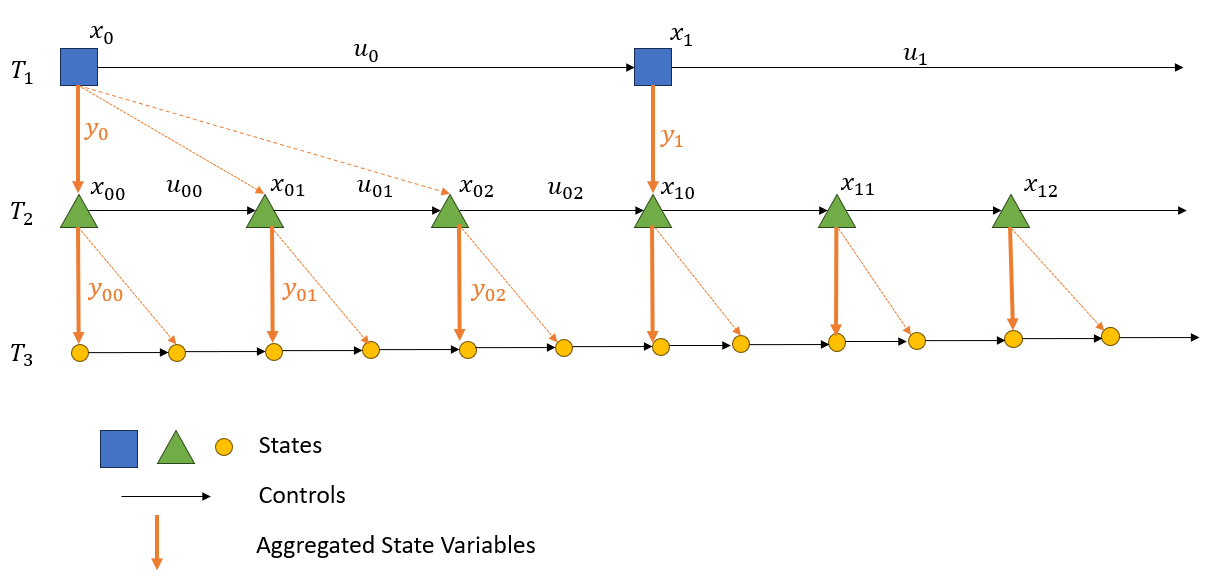}
    \caption{Connection between decision variables across different timescales.}
    \label{fig:multiscale_vis}
\end{figure}

A visualization of the interaction of variables can be found in Figure \ref{fig:multiscale_vis}.  In the figure, each node is a random decision variable. Rows depict within-layer state dynamics, where aggregated state variables are drawn between the rows.

\subsection{Non-anticipativity under Multi-timescale Settings} \label{section:nac_multitimescale}

We conclude the definition of a multi-timescale problem by defining the concept of non-anticipativity in a multi-timescale setup. 

Non-anticipativity means that the decisions made at each tick do not use future information. In this framework, decisions are made in the dictionary order of tick indices. Consider a three-timescale problem; for example, the order is: \((0), (0,0), (0,0,0), (0,0,1), (0,0,2), \dots, (0,1), (0,1,0), \dots\). When we are ready to make a decision at tick \((0,0,2)\), we have access to the random variables \(\omega\) with indices \((0), (0,0), (0,0,0), (0,0,1), (0,0,2)\), which captures all information available up to that point.

However, there are situations where ticks correspond to the same physical time, e.g., ticks \((0), (0,0), (0,0,0)\) all represent 0:00 AM in our motivational example. This can lead to different interpretations of non-anticipativity:

\begin{itemize}
    \item \textbf{Interpretation 1:} When several decisions are made at the same physical time and no new information is revealed between these decisions, they must satisfy the same non-anticipativity constraints.
    \item \textbf{Interpretation 2:} Although the decisions occur simultaneously, decisions across different timescales are made sequentially. For instance, after a decision is made at a slower timescale, new information may be revealed that can be used to inform subsequent decisions at a faster timescale.
\end{itemize}

To formally capture the notion of ``the available information so far,'' we define the set of all ticks by \(\mathcal{T} = \mathcal{T}_1 \cup \mathcal{T}_2 \cup \mathcal{T}_3\). By sorting the tuples in dictionary order, we establish a total order \(\leq\) on \(\mathcal{T}\). The difference between the two interpretations is in the treatment of ticks that represent the same physical time. In the first interpretation, such ticks are considered ``equal" to each other, implying they share the same information and constraints. In contrast, in the second interpretation, these ticks are not considered equal. That means even at the same physical time, the decisions are still made sequentially, and new information can be injected in between and influence decisions.

Regardless of the interpretation chosen, we can define a natural filtration:

\begin{equation}
    F_{\tau} = \sigma(\{\omega_{\tau'}: \tau' \leq \tau\})
\end{equation}

In this definition, the index set \(\{\tau': \tau' \leq \tau\}\) refers to the set of all ticks up to (and including) the tick \(\tau\), similar to the integer ordering of stage numbers in the multistage stochastic programming setting in Section \ref{section:nac_single_timescale}. Non-anticipativity in the multi-timescale setting means $x_\tau$, $u_\tau$ are measurable with respect to $F_{\tau}$ for each tick $\tau \in \mathcal T$.

We clarify that in our motivational example and all subsequent descriptions of the formulations, we use the first interpretation. We assume that ticks representing the same physical time are considered ``equal,'' meaning they share the same information and are subject to the same non-anticipativity constraints.

\section{Equivalent Multistage Problem of a Multi-timescale Problem} \label{section:equivalent_multistage_problem}

\begin{figure}
        \centering
        \includegraphics[width=\linewidth]{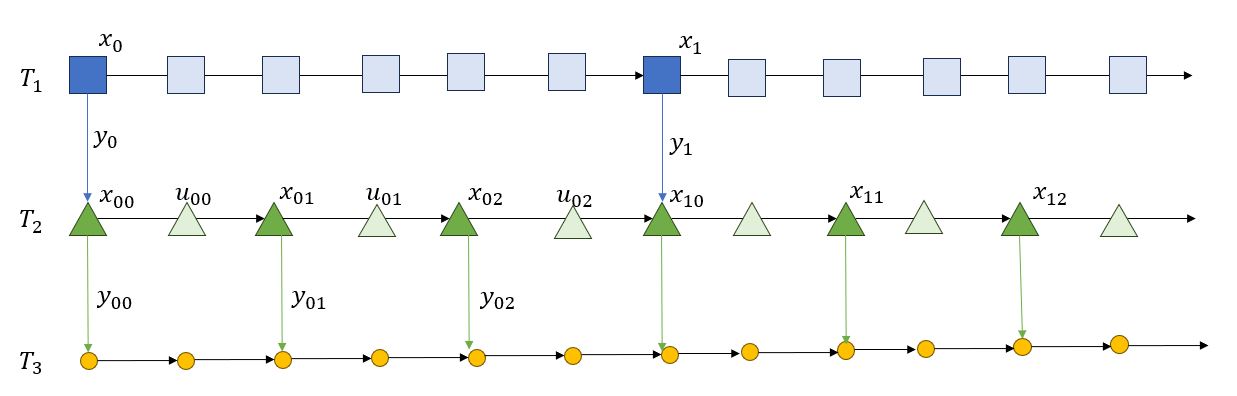}
        \caption{Visualization of equivalent multistage formulation (random variable view)}
        \label{fig:equivalent_multistage_full}
\end{figure}

With the non-anticipativity definition, our multi-timescale problem can be reformulated as an equivalent multistage problem. A representation of this is provided in Figure \ref{fig:equivalent_multistage_full}. In this section, we discuss why this multi-timescale structure creates challenges for traditional multistage algorithms.

\begin{figure}
    \centering
    \includegraphics[width=\linewidth]{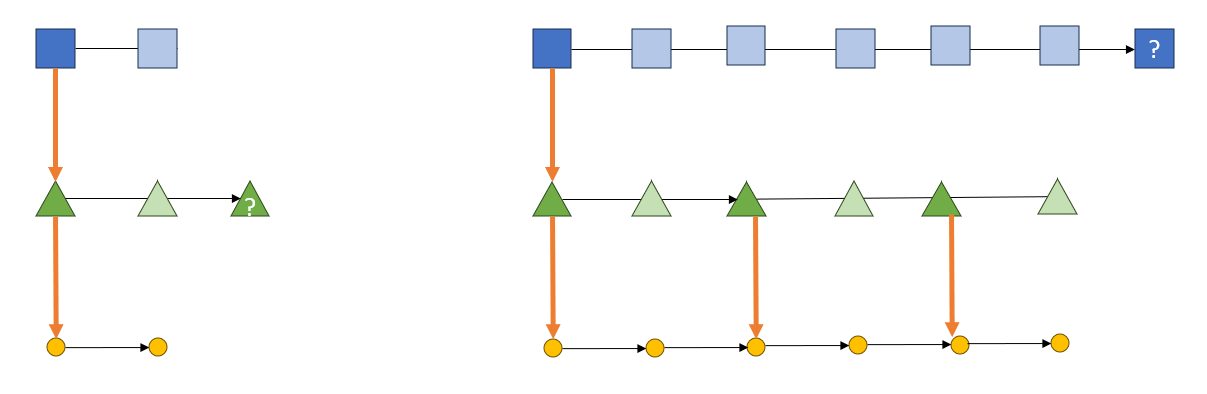}
    \caption{Intermediate stage in equivalent multistage formulation (random variable view)}
    \label{fig:equivalent_multistage_intermediate}
\end{figure}

Observe that the resulting multistage problem cannot be directly formulated as a traditional optimal control problem. This is because of the unique dependency structure within the multi-timescale framework: the state in a given timescale may depend on a previous tick within the same timescale. However, when translated to the equivalent multistage problem, these dependencies could span multiple stages backward, crossing multiple timescales. To further elaborate, we visualize two intermediate stages in this equivalent multistage problem in Figure \ref{fig:equivalent_multistage_intermediate}. In the left subfigure of Figure \ref{fig:equivalent_multistage_intermediate}, the dynamics of the state node \(x_{01}\) (marked as a green triangle with a "?") are influenced by \(x_{00}\) and the aggregated state variables \(y_0\). Notably, this decision is made only after the faster ticks \(x_{000}\) and \(x_{001}\) that occur in between. On the right side of Figure \ref{fig:equivalent_multistage_intermediate}, the situation is even more extreme; the node \(x_1\) (marked with a blue square and "?") requires information about states from several periods prior.

If we decide to solve this as a multistage problem using resource-directive decomposition, a potential remedy would be to expand the definition of state variables to encompass states across all timescales and modify the state dynamics accordingly. However, this approach leads to a significant increase in the number of stages and unnecessarily expands the dimension of the state variables. It is well-documented that the tractability of resource-directive decompositions relies on maintaining a small dimension of state variables, as large dimensions can result in the "curse of dimensionality." This expansion would likely limit the effectiveness of the decomposition strategy.

\begin{figure}
    \centering
    \includegraphics[width=0.8\linewidth]{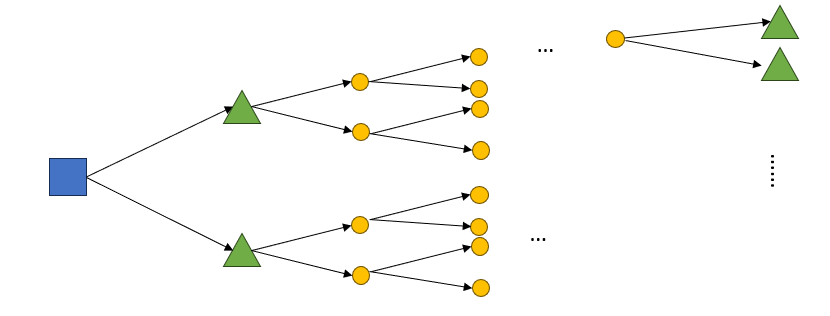}
    \caption{Equivalent multistage formulation in scenario tree view}
    \label{fig:equivalent_multistage_scenario_tree}
\end{figure}

This multistage problem also poses a challenge for (scenario-based) price-directive decomposition. This is due to the number of scenarios required in a scenario tree growing exponentially with the number of stages. Even a modest multi-timescale problem can result in a huge scenario tree, as shown in Figure \ref{fig:equivalent_multistage_scenario_tree}.  Consequently, scenario-based methods likely cannot solve this multistage equivalent problem beyond relatively small, illustrative examples.

\section{Improve Tractability by Approximation} \label{section:tractability}
Our approximation is inspired by techniques proposed by Kaut \cite{kaut_multi-horizon_2014} and Maggioni \cite{maggioni_bounds_2020}. In this section, we first review their contributions and why their method cannot be directly used in our motivational example. Subsequently, we introduce our approximation method, termed the ``synchronized state approximation,'' followed by some discussion on this method.

\subsection{Multihorizon Scenario Tree \cite{kaut_multi-horizon_2014}}
\begin{figure}
    \centering
    \includegraphics[width=0.6\linewidth]{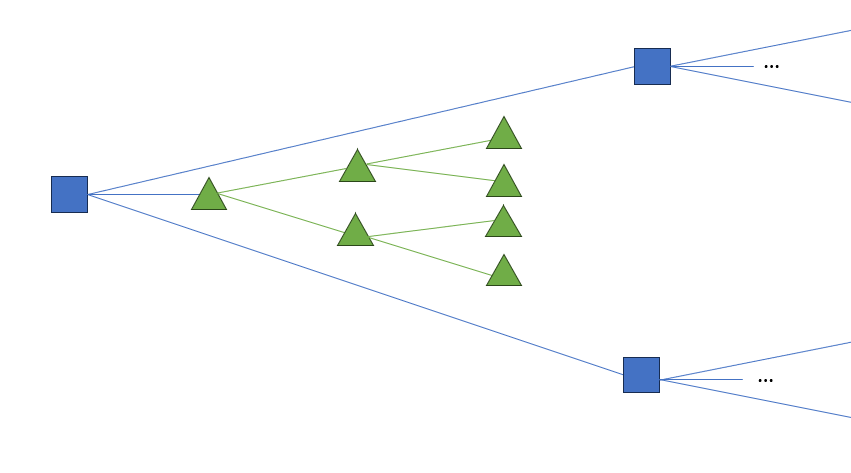}
    \caption{Kaut's Multihorizon Scenario Tree \cite{kaut_multi-horizon_2014} (in scenario view)}
    \label{fig:kaut_mh_tree}
\end{figure}

 In Kaut's work, the following assumptions are needed:
\begin{enumerate}
    \item Uncertainty in the fast timescale is independent of uncertainty in the slow timescale.
    \item Decisions in a slow timescale do not affect decisions in a fast timescale.
    \item The first decision in the fast timescale cannot depend on the previous decisions in the fast timescale.
\end{enumerate}
If these assumptions are met, a simplified scenario tree can be constructed. Figure \ref{fig:kaut_mh_tree} illustrates a two-timescale example of this simplified tree. In this figure, the fast nodes (green triangle-shaped \img{triangle.PNG}) are depicted as directly ``hanging'' from the slow nodes (blue square-shaped \img{square.PNG}), rather than branching independently at each occurrence as shown in Figure \ref{fig:equivalent_multistage_scenario_tree}. This significantly improves tractability because the total number of scenarios required to construct this tree is the product of the number of scenarios needed for each timescale.

In our approach, Assumption 1 can be somewhat relaxed by using scenario descriptions at the slowest timescale, which allows for multiple lag dependencies. Assumption 2 is addressed by defining aggregated state variables, which communicate information from the slow to the fast timescale. However, Assumption 3 is the key reason why the tree can be simplified but cannot be justified in our motivational example (refer to Figure \ref{fig:synchronized_states_assumption} and the subsequent discussion). Consequently, we propose an approximation termed the ``synchronized states approximation'' in the next subsection.

\subsection{Synchronized State Approximation}

\begin{figure}
    \centering
    \includegraphics[width=\linewidth]{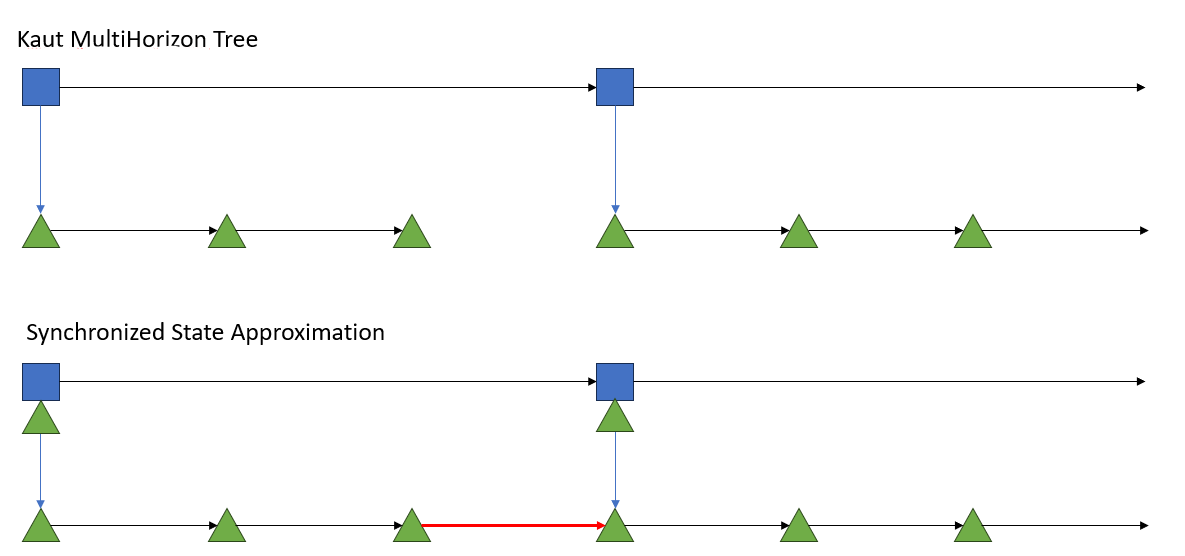}
    \caption{Assumption in Kaut's Multihorizon Tree in random variable view and synchronized states approximation}
    \label{fig:synchronized_states_assumption}
\end{figure}

Assumption 3 in Kaut's work significantly simplifies the handling of state variables by removing the need to maintain a coherent state across the entire time horizon. Instead, coherence is only required within each problem segment — specifically, between the fast ticks that correspond to the same parent tick. The upper subfigure of Figure \ref{fig:synchronized_states_assumption} shows this dependency. In the figure, there is no connection between the third green triangle and the fourth green triangle. In the optimal solution of the equivalent multistage problem, the slow states, indicated by the second blue square in the first row, remain the same across all scenarios. A detailed and robust treatment can be found in the referenced paper \cite{kaut_multi-horizon_2014}.

However, this assumption is not directly applicable to our motivational example because state dynamics need to be kept across the entire time horizon. This consistency issue can be resolved by making sure that states in the fast timescale are made to be equal at the beginning of each time segment; this is called ``synchronized states''. One practical approximation to address this issue is to artificially align the system states in the fast timescale at segment starts by modifying the non-anticipativity constraints.

Equivalently, this approximation involves duplicating the state from a fast tick to match the timescale of a slow tick, as depicted in the lower part of Figure \ref{fig:synchronized_states_assumption}, where a green triangle is directly aligned beneath a blue square. At the end of each segment, a feasible control, indicated by an orange arrow, must be implemented to adjust the state so that it matches the prescribed conditions at the slow timescale.

Although this clears the consistency issues of the state dynamics, this approximation imposes artificial constraints on the fast timescales. The operational flexibility of the systems is reduced in the fast timescale. Therefore, assuming a minimization problem, the added constraint likely increases the objective value. 

Another limitation is that the resulting synchronized states approximation problem may be infeasible. Therefore, this method may not be suitable for some applications.

\section{Instantiation} \label{section:instantiation}

In the previous sections, we described how to model a multi-timescale problem using random decision variables. However, to actually solve the problem, we need to transform it into a form where these decision variables become standard optimization variables. This section discusses strategies for making that transformation.

\subsection{Scenario-based Instantiation} \label{section:scenario_based_instantiation}

\begin{figure}
    \centering
    \includegraphics[width=0.5\linewidth]{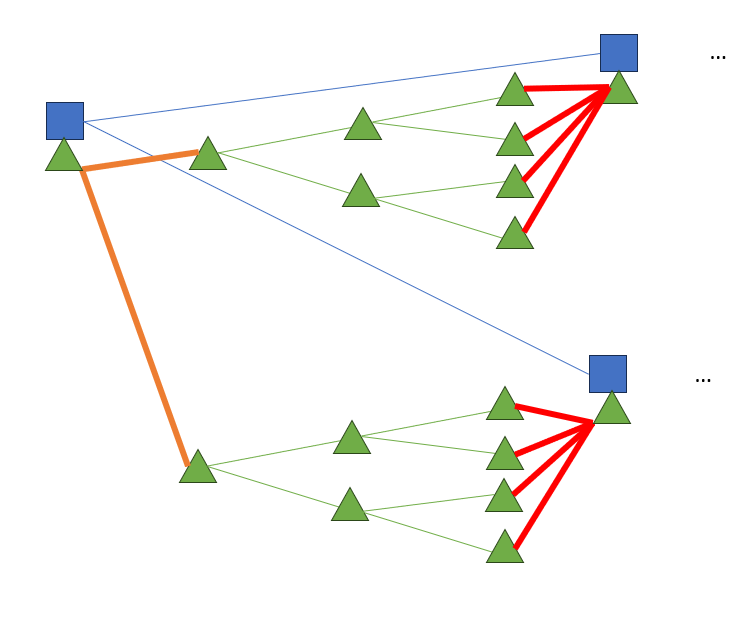}
    \caption{A two-timescale example with a scenario-based instantiation. The first timescale node is marked with blue squares, and the green triangles mark the nodes of the second timescale. There are two periods in the first timescale and two periods in the second timescale. }
    \label{fig:multiscale_example}
\end{figure}

A direct way of instantiating the problem is to form scenario trees in each timescale and pose the problem in a similar way as in the multi-horizon trees \cite{maggioni_bounds_2020}.

Here, we use a two-timescale example for illustration. we let $T_1 = (V_1, E_1)$ be a rooted tree representing the first timescale scenario tree. The probability of going through a node is notated as $p_{n_1}$ for $n_1 \in V_1$. Under the synchronization approximation at each node $n_1 \in V_1$ we need to define these variables: $x_{n_1}$, the state in timescale 1; $u_{n_1}$, the control taken in timescale 1; $y_{n_1}$, the aggregate state variable generated in timescale 1; $z_{n_1}$, the prescribed state for timescale 2.

To complete the second timescale, between each of the two first timescale scenario nodes $(m_1, n_1) \in E_1$, we define another (rooted) scenario tree $T_2(n_1) = (V_2(n_1), E_2(n_1))$. We should be aware that the depth of the depth of the second timescale scenario tree has one more depth than the number of random variables (see example in Figure \ref{fig:multiscale_example}.) Depending on the formulation (decision-hazard, or hazard-decision), the root node may be deterministic.

Similarly we have the probability going through each node $p(n_2)$ for $n_2 \in V_2(n_1)$, and define variables $x_{n_2}, u_{n_2}$. Note since there are no finer timescales, we do not need $y_{n_2}$ and $z_{n_2}$.

\begin{align}
\min & \sum_{n_1 \in N_1} p_{n_1} f_{n_1}(x_{n_1}, u_{n_1}) + \sum_{(m_1, n_1) \in E_1} \sum_{n_2 \in V_2(n_1)} p_{n_2} f_{n_2}(x_{n_2}, u_{n_2}) \\
& \textbf{Timescale 1 Formulation} \nonumber \\ 
& x_{n_1} = D_{n_1}(x_{m_1}, u_{m_1}) \qquad \forall (m_1, n_1) \in E_1 \\
& u_{n_1} \in U_{n_1}(x_{m_1}) \qquad \forall (m_1, n_1) \in E_1\\
& y_{n_1} = C_{n_1}(x_{n_1}) \qquad \forall n_1 \in V_1 \\
& \textbf{Timescale 2 Formulation} \nonumber \\ 
& x_{n_2} = D_{n_2}(x_{m_2}, u_{n_2}, y_{m_1}) \qquad \forall (m_1, n_1) \in E_1 \forall (m_2, n_2) \in E_2(n_1) \\
& u_{n_2} \in U_{n_2}(x_{m_2}, y_{m_1}) \qquad \forall (m_1, n_1) \in E_1 \forall (m_2, n_2) \in E_2(n_1) \\
& x_{n_2} = z_{m_1} \qquad \forall (m_1, n_1) \in E_1 \text{ and } n_2 \text{ is root of } T_2(n_1) \label{eqn:example_start_state} \\
& x_{n_2} = z_{n_1} \qquad \forall (m_1, n_1) \in E_1 \text{ and } n_2 \text{ is a leaf of } T_2(n_1) \label{eqn:example_ending_state}
\end{align}

The objective of the overall problem is the sum of the objectives in these two scales. It is worthy to note that even $z_{m_1}$ does not appear in the objective of the first timescale, they belong to that timescale, because the second timescale will need to use those shared variables.

Using this technique, we achieve the same reduction in problem size as the multihorizon tree approach. Consider a three-timescale problem: the first timescale contains $T_1$ ticks, each tick in the first timescale corresponds to $T_2$ ticks in the second timescale, and each tick in the second timescale corresponds to $T_3$ ticks in the third timescale. To represent the tree structure for the first timescale, approximately $e^{T_1}$ nodes are required. Similarly, constructing the tree for the second timescale requires $e^{T_2}$ nodes for each node in the first timescale, and constructing the third timescale tree requires $e^{T_3}$ nodes for each node in the second timescale. Consequently, the entire multihorizon scenario tree would require roughly $e^{T_1 + T_2 + T_3}$ nodes. We see multi-horizon trees need significantly number of nodes compared to the expanded multistage problem in Section \ref{section:equivalent_multistage_problem}, which would need about $e^{T_1 T_2 T_3}$ nodes. This is a meaningful step toward tractability.

However, even with the multihorizon scenario tree, the number of optimization variables and constraints is still huge. For instance, in a small model comprising one slow-ramping generator, one fast-ramping generator, and one renewable generator with heavily simplified dynamics, 714,177 variables and 1,306,045 constraints were needed. Hence full scenario-based instantiation is not a good strategy for large-scale problems.

\subsection{Specialized Value Functions} \label{section:value_function_instantiation}


To address this scalability challenge outlined above, we propose an alternative instantiation based on the principles of dynamic programming (DP). Our formulation extends the traditional DP framework: where standard value functions depend only on the current state, ours must also account for both the prescribed target state at the end of the segment and the fixed aggregated state passed down from the parent timescale. This approach is effective for timescales where the underlying uncertainty follows a Markovian structure, allowing for a more compact problem representation. To construct this DP formulation, we must first decompose the problem into a sequence of smaller, self-contained segment subproblems.

A \textbf{segment} is a consecutive sequence of ticks within a given timescale, say timescale $j$, that is entirely governed by a single tick from its immediate parent timescale, $j-1$. Specifically, for any parent tick $\tau_{j-1}$, the corresponding segment in timescale $j$, denoted $S(\tau_{j-1})$, comprises all ticks from the starting tick $\mathcal{T}_j(\tau_{j-1})$ up to, but not including, the next starting tick $\mathcal{T}_j(\tau_{j-1}+)$. A critical feature of a segment is that all decisions within it are influenced by the same aggregated state variable, $y_{\tau_{j-1}}$, passed down from the parent tick. This structure allows us to view the optimization over a segment as a self-contained subproblem, which must satisfy boundary conditions imposed by the slower timescale.







With the concept of a segment established, we can now formally define the value function for our subproblem. This function, denoted $Q_{\tau_j}(x_{\tau_j}, \hat x_{\mathcal T_j(\tau_{j-1}+)}, y_{\tau_{j-1}})$, represents the minimum expected cost to optimally steer the system from the current tick $\tau_j$ through the remainder of the segment. A key aspect of this function is to enforce a seamless ``hand-off'' between segments. This is achieved by enforcing a boundary condition: the state resulting from the segment's final action must equal the prescribed synchronized \textbf{goal state}, $\hat x_{\mathcal T_j(\tau_{j-1}+)}$, which is the required initial state for the next segment.


Mathematically, the cost-to-go is defined via recursion. The terminal condition at the end of the segment ensures that the ending state matches the target ending state. It is expressed as:

\begin{equation} \label{eqn:ending_condition}
Q_{\tau_j}(x_{\tau_j}, \hat x_{\mathcal T_j(\tau_{j-1}+)} , y_{\tau_{j-1}}) = \mathbf{1}\{\hat x_{\tau_j} = x_{\mathcal T_j(\tau_{j-1}+)}\}
\qquad \text{if} \qquad \tau_j = \mathcal T_j(\tau_{j-1}+)
\end{equation}

Alternatively, L1 penalty $\rho_{\tau_j} \| x - \tilde x\|$ can be applied when the process fails to drive the state to the prescribed ending state. This can be viewed as a relaxation of \eqref{eqn:ending_condition}.

For all other ticks in the segment, the value function is defined recursively using the Bellman principle:

\begin{align} 
Q_{\tau_j}(x_{\tau_j}, \hat x_{\mathcal T_j(\tau_{j-1}+)} , y_{\tau_{j-1}}) = \min & \ f_{\tau_j}(x_{\tau_j}, u_{\tau_j}) \label{eqn:value_function1} \\
& + Q_{\tau_{j+1}}^{\tau_{j}}(x_{\mathcal T_{j+1}(\tau_j)}, \hat x_{\mathcal T_{j+1}(\tau_j +)}, y_{\tau_{j}}) \label{eqn:value_function2} \\
& + E[Q_{\tau_j}(x_{\tau_j+}, x_{\mathcal T_j(\tau_{j-1}+)} , y_{\tau_{j-1}})]    \label{eqn:value_function3} 
\end{align}

\[
\text{ subject to 
\eqref{eqn:agg_state_generation},\eqref{eqn:multiscale_feasible_actions},\eqref{eqn:multiscale_transition}}
\]

In the definition above, the terms numbered \eqref{eqn:value_function1} and \eqref{eqn:value_function3} represent the usual value functions. The term \eqref{eqn:value_function2} recursively includes the value function from the next, more granular timescale.

This instantiation approach exploits the independence structure to achieve a more compact representation of the problem. Approximate DP-based methods, such as SDDP \cite{powell_approximate_2011}, are well-suited to the resulting problem. Notably, this approach allows flexibility: we can combine different instantiation methods across timescales, shown through the example in the next section.

\section{A Three-Timescale Stochastic Unit Commitment Example}
\label{sec:multi-timescale_instance}

To demonstrate the use of different instantiation methods, this section construct a concrete mathematical model for the stochastic unit commitment problem first discussed as a motivational example in Section~\ref{sec:motivational_example}.

In this motivational problem, scenario-based instantiation should be applied to the first (DA) timescale to capture autocorrelation, while the second (ST) and third (RT) timescales can use value function instantiation to exploit the Markovian property of the randomness. In applications employing stochastic bridge-based modeling \cite{glanzer_multiscale_2020}, value function-based instantiation is a more natural choice.

It is important to distinguish between two key temporal concepts used in the actual modeling:
\begin{enumerate}
\item \textbf{Decision Resolution}: This is the time interval between consecutive control actions. For example, a one-hour decision resolution means that operators can make adjustments every hour.
\item \textbf{Stochastic Stage Length}: This defines how frequently new information about uncertain variables is revealed. In a scenario tree, this corresponds to the time between branching points.
\end{enumerate}

The principle of non-anticipativity requires that decisions be based only on information that has already been revealed. In our Day-Ahead (DA) formulation, we use a 1-hour decision resolution but a 4-hour stochastic stage length. This means that at the start of each 4-hour block, the model determines a sequence of four hourly decisions based only on the information known up to that point.

\textbf{Important Remark:} This structure implies that the model has \textbf{perfect foresight of the day-ahead information ($\omega^{(1)}$) within each 4-hour segment}. In our application, this means the DA model knows the evolution of its coarse renewable forecast for the next four hours. However, it \textbf{cannot see the more refined forecast information} that is only revealed to the faster, short-term (ST) timescale. The DA model must therefore collaborate with the ST timescale—via the value function $V^{(2)}$—to make decisions with this limited foresight, anticipating that the ST model will handle the finer-grained uncertainty.

While reducing the stochastic stage length to one hour would remove this specific look-ahead assumption, it would also require the scenario tree to branch much more frequently, severely impacting computational tractability.

\begin{table}[htbp]
\centering
\caption{Characteristics of the Three Timescales in the Unit Commitment Model}
\label{tab:timescale_characteristics}
\begin{tabular}{@{}lllll@{}}
\toprule
\textbf{Timescale} & \textbf{Primary Decisions} & \textbf{Decision Resolution} & \textbf{Stochastic Stage Length} & \textbf{Uncertainty Model} \\ 
\midrule
DA & Slow-ramping units & 1 hour & 4 hours & Scenario Tree \\
ST & Fast-ramping units & 1 hour & 1 hour & Scenario Lattice \\
RT & Economic dispatch & 15 min & 15 min & Independent Noise \\ 
\bottomrule
\end{tabular}
\end{table}

We use different instantiation methods for each timescale to balance modeling fidelity with computational tractability.

For the DA timescale, we employ the scenario-based approach from Section~\ref{section:scenario_based_instantiation}. This allows for a detailed representation of the slow-ramping generators, including history-dependent constraints, notably minimum uptime and downtime.

For the Short-Term (ST) and Real-Time (RT) timescales, we use the value function approach from Section~\ref{section:value_function_instantiation} to accommodate the scenario lattice structure. This choice necessitates a simplification for the fast-ramping generators in the ST timescale. To use value functions, the system dynamics must be Markovian (i.e., depend only on the current state). Therefore, we omit history-dependent constraints for these generators. This simplification is physically justified, as fast-ramping generators typically have very short or negligible minimum uptime and downtime requirements; generators with significant history-dependent constraints would be more appropriately modeled in the DA timescale. From a computational perspective, the alternative of expanding the state space to track their history would be prohibitive.

Finally, in the RT timescale, renewable sources are treated as ``profiled generators,'' a method also used in the \texttt{UnitCommitment.jl} package \cite{santos_xavier_unitcommitmentjl_2024}. In this approach, we only model their power output, which is constrained by generation limits provided as input data.

The following tables provide a high-level overview of the modeling components applied to each timescale.

\begin{threeparttable}
    \caption{Overview of Generator Modeling Components by Timescale (Operational Constraints)}
    \label{tab:modeling_overview_op}
    \begin{tabular}{@{}l l c c c c c@{}}
        \toprule
        \textbf{Timescale} & \textbf{Generators} & \textbf{On/off} & \textbf{Power} & \textbf{Gen. Limits} & \textbf{Ramping Limit} & \textbf{Min Up/Down} \\
        \midrule
        DA & Slow-ramping & Y & Y & Y & Y & Y \\
        ST & Fast-ramping & Y & Y & Y & Y & N \\
        RT & Profiled & N/A & Y & Y\tnote{a} & N/A & N/A \\
        \bottomrule
    \end{tabular}
    \begin{tablenotes}
        \item[a] Generation limits are not decision variables but are provided as exogenous data / uncertainty.
    \end{tablenotes}
\end{threeparttable}

\begin{threeparttable}
    \caption{Overview of Generator Modeling Components by Timescale (Cost Components)}
    \label{tab:modeling_overview_cost}
    \begin{tabular}{@{}l l c c c@{}}
        \toprule
        \textbf{Timescale} & \textbf{Generators} & \textbf{Prod. Cost} & \textbf{No-load Cost} & \textbf{Start/Shut Cost} \\
        \midrule
        DA & Slow-ramping & Y & Y & Y \\
        ST & Fast-ramping & Y & Y & Y\tnote{a} \\
        RT & Profiled & N/A & N/A & N/A \\
        \bottomrule
    \end{tabular}
    \begin{tablenotes}
        \item[a] Start-up cost is a fixed value that does not account for warm-start conditions.
    \end{tablenotes}
\end{threeparttable}

\subsection{Day Ahead Timescale}

\begin{table}[htbp]
    \centering
    \caption{Notation: Sets}
    \label{tab:notation_sets}
    \begin{tabular}{@{}ll@{}}
        \toprule
        \textbf{Set} & \textbf{Description} \\
        \midrule
        $\mathcal G^{(1)}, \mathcal G^{(2)}, \mathcal G^{(3)}$ & Sets of generators in timescale 1 (slow-ramping), 2 (fast-ramping), and 3 (real-time/renewable). \\
        $\mathcal B$ & Set of buses. \\
        $\mathcal G^{(1)}(b), \mathcal G^{(2)}(b), \mathcal G^{(3)}(b)$ & Sets of generators of type 1, 2, and 3 at bus $b \in \mathcal B$. \\
        $\mathcal T^{(1)}, \mathcal T^{(2)}, \mathcal T^{(3)}$ & Sets of ticks for timescale 1 (4-hour), 2 (1-hour), and 3 (15-min). \\
        $\mathcal S^{(1)}$ & Set of scenarios in timescale 1. \\
        \bottomrule
    \end{tabular}
\end{table}

\begin{table}[htbp]
    \centering
    \caption{Notation: Timescale 1 (DA) Variables and Parameters}
    \label{tab:notation_timescale1}
    \begin{tabular}{@{}p{0.2\textwidth} p{0.7\textwidth}@{}}
        \toprule
        \textbf{Symbol} & \textbf{Description} \\
        \midrule
        $\pi^s$ & Probability of scenario $s$ in timescale 1. \\
        $\omega^{(1),s}$ & A specific scenario; one possible realization of the random outcomes $(\omega^{(1),s}_{\tau})_{\tau \in \mathcal T^{(1)}}$ over the planning horizon. \\

        \addlinespace
        $\mathbf{x}^{(1),s}$ & Vector of state variables for scenario $s$. For each 4-hour tick $\tau \in \mathcal T^{(1)}$, it contains hourly values for: \\
        \quad $x_{g,\tau}^s$ & On/off status (binary) of \textbf{slow-ramping} generator $g \in \mathcal G^{(1)}$. \\
        \quad $p_{g,\tau}^s$ & Power output of \textbf{slow-ramping} generator $g \in \mathcal G^{(1)}$. \\
        \addlinespace
        $\mathbf{y}^{(1),s}$ & Aggregate state variable for scenario $s$. \\
        \quad $y_{b, \tau}^{(1),s}$ &  Total power committed at each bus $b$ for each tick $\tau$. \\
        \addlinespace
        $\mathbf{u}^{(1),s}$ & Vector of control variables for scenario $s$. For each tick $\tau \in \mathcal T^{(
        1)}$, it contains: \\
        \quad $v_{g,\tau}^s$ & Startup decision for generator $g$. \\
        \quad $w_{g,\tau}^s$ & Shutdown decision for generator $g$. \\
        \quad $\Delta p_{g,\tau}^s$ & Ramping decision for generator $g$. \\
        \addlinespace
        $\mathcal X^{(1)}$ & Generator polytope for all slow-ramping generators; the set of mixed-integer constraints defining the feasible operating region for each generator $g \in \mathcal G^{(1)}$. \\
        \addlinespace
        $\mathbf{z}^{(1),s}$ & \textbf{Prescribed target state} for fast ramping generators. \\
        \quad $x_{g,\tau}^s$ & \textbf{Prescribed} on/off status (binary) of \textbf{fast-ramping} generator $g \in \mathcal G^{(2)}$. \\
        \quad $p_{g,\tau}^s$ & \textbf{Prescribed} power output of \textbf{fast-ramping} generator $g \in \mathcal G^{(2)}$. \\
        \bottomrule
    \end{tabular}
\end{table}

The day ahead timescale minimization problem can be defined as follows.

\begin{align}
\min_{\mathbf{x,u,y,z}} \quad & \sum_{s \in \mathcal S^{(1)}}\pi^{s} \left( \mathbf{c}^{(1),s} \mathbf{x}^{(1),s} + \mathbf{d}^{(1),s} \mathbf{u}^{(1),s} +\sum_{\tau \in \mathcal T^{(1)}} V^{(2)}_\tau(\mathbf{z}^{(1),s}_\tau,\mathbf{z}^{(1),s}_{\tau+}, \mathbf{y}^{(1),s}_\tau, \omega^{(1),s}_{\tau}) \right) \label{eqn:DA_obj}\\
\text{s.t.} \quad & (\mathbf{x}^{(1),s}, \mathbf{u}^{(1),s}) \in \mathcal X^{(1)}, \quad \forall s \in \mathcal S^{(1)}, \label{eqn:DA_feasible_trajectory} \\
& y^{(1),s}_{b,\tau} = \sum_{g \in \mathcal G^{(1)}(b)} p_{g,\tau}^{s}, \quad \forall s \in \mathcal S^{(1)}, b \in \mathcal B, \tau \in \mathcal T^{(1)}, \label{eqn:DA_agg_state} \\
& \mathbf{z}^{(1),s}_\tau \in X^{(2)}, \quad \forall s \in \mathcal S^{(1)}, \tau \in \mathcal T^{(1)}, \label{eqn:DA_feasible_prescribed_ST_states}\\
& \mathbf{\bar x}^{(1)} - N^{(1)} \mathbf{x}^{(1)} = 0, \label{eqn:DA_NAC1} \\
& \mathbf{\bar z}^{(1)} - N^{(1)} \mathbf{z}^{(1)} = 0. \label{eqn:DA_NAC2}
\end{align}

The terms $\mathbf{c}^{(1),s} \mathbf{x}^{(1),s} + \mathbf{d}^{(1),s} \mathbf{u}^{(1),s}$ in the objective \eqref{eqn:DA_obj} represent the operational costs of the slow-ramping generators. These costs, combined with the polytope constraint $(\mathbf{x}^{(1),s}, \mathbf{u}^{(1),s}) \in \mathcal X^{(1)}$ from \eqref{eqn:DA_feasible_trajectory}, define a unit commitment subproblem for the timescale-1 system in scenario $s$. This polytope $\mathcal X^{(1)}$ encapsulates all operational constraints. Because this problem is formulated on a per-scenario basis, $\mathcal X^{(1)}$ can readily model all types of history-dependent relationships, most notably the minimum uptime and downtime constraints that link decisions across multiple time periods. For practical implementation, the polytope can be constructed using specialized packages that generate computationally advantageous tight formulations, which are effective because they leverage these links between variables across different ticks.

The term $V^{(2)}_\tau$ in the objective function \eqref{eqn:DA_obj} represents the future costs incurred in the faster, short-term (ST) timescale. Specifically, it is a value function that provides the optimal expected cost to operate the timescale-2 system over the segment from day-ahead tick $\tau$ to $\tau+$. This function is parameterized by day-ahead decisions and parameters:

\begin{itemize}
    \item The prescribed states $\mathbf{z}^{(1),s}_\tau$ and $\mathbf{z}^{(1),s}_{\tau+}$ act as boundary conditions, defining the required state of the fast-ramping generators at the start and end of the segment.
    \item The ST subproblem operates using the aggregated power $\mathbf{y}^{(1),s}_\tau$ made available from timescale 1 and is guided by the coarse renewable generation forecast $\omega^{(1),s}_{\tau}$ provided at the day-ahead stage.
    \item To ensure the problem is well-posed, constraint \eqref{eqn:DA_feasible_prescribed_ST_states} restricts these prescribed states to be physically feasible within the operational limits of the timescale-2 generators.
\end{itemize}

Finally, constraints \eqref{eqn:DA_NAC1} and \eqref{eqn:DA_NAC2} enforce non-anticipativity on the state variables ($\mathbf{x}^{(1)}$) and the prescribed states ($\mathbf{z}^{(1)}$). In this formulation, $\mathbf{\bar x}^{(1)}$ and $\mathbf{\bar z}^{(1)}$ represent the single, implementable policy decisions. The constraints link the individual scenario-based decisions to these implementable variables, forcing all scenarios that are indistinguishable at a given 4-hour stochastic stage to have identical decisions.

\subsection{Short Term Timescale}

\begin{table}[htbp]
    \centering
    \caption{Notation: Timescale 2 (ST) Variables}
    \label{tab:notation_timescale2}
    \begin{tabular}{@{}p{0.2\textwidth} p{0.7\textwidth}@{}}
        \toprule
        \textbf{Symbol} & \textbf{Description} \\
        \midrule
        $\mathbf{x}^{(2)}_\tau$ & Vector of state variables at tick $\tau \in \mathcal T^{(2)}$. \\
        & \quad $x_{g,\tau}$: On/off status for generator $g \in \mathcal G^{(2)}$. \\
        & \quad $p_{g,\tau}$: Power output for generator $g \in \mathcal G^{(2)}$. \\
        \addlinespace
        $\mathbf{u}^{(2)}_\tau$ & Vector of control variables at tick $\tau \in \mathcal T^{(2)}$. \\
        & \quad $v_{g,\tau}$: Startup decision for generator $g$. \\
        & \quad $w_{g,\tau}$: Shutdown decision for generator $g$. \\
        & \quad $\Delta p_{g,\tau}$: Ramping decision for generator $g$. \\
        \addlinespace
        $\mathbf{y}^{(2)}_\tau$ & Aggregate state variable. \\
        & \quad $y^{(2)}_{b,\tau}$Total power committed at each bus $b$, including both slow-ramping and fast-ramping generators. \\
        \addlinespace
        $X^{(2)}$ & The set of feasible states (e.g., generation limits) for all fast-ramping generators $g \in \mathcal G^{(2)}$. \\
        \bottomrule
    \end{tabular}
\end{table}

The short-term (ST) timescale problem is defined by a value function $V^{(2)}_\tau$. To avoid confusion, we first define an auxillary process, $Q^{(2)}_\tau$ via the following dynamic programming recursion, which has another parameter, attached to the end, to track the Markovian random outcome via the appended $\omega^{(2)}_{\tau}$:

If $\tau \in \mathcal T^{(2)}$ is not at the end of the segment:
\begin{align}
    Q^{(2)}_\tau(\mathbf{x}_\tau, \mathbf{\tilde x}, \mathbf{y}^{(1)};\omega^{(1)}_{\tau_1}, \omega^{(2)}_{\tau}) = \min \ &\mathbf{c}^{(2)}_\tau \mathbf{x}_\tau  + \mathbf{d}^{(2)}_\tau \mathbf{u}_\tau \label{eqn:ST_immediate_cost}  \\
    &+ E_{\omega^{(3)}_\tau | \omega^{(2)}_\tau}[V^{(3)}_\tau(\mathbf{y}_\tau^{(2)}, \omega^{(2)}_\tau)] \label{eqn:ST_RT_cost} \\
    &+ E_{\omega^{(2)}_{\tau+} | \omega^{(1)}_{\tau_1},\omega^{(2)}_\tau}[Q^{(2)}_{\tau+}(\mathbf{x}_{\tau+}, \mathbf{\tilde x}, \mathbf{y}^{(1)};\omega^{(1)}_{\tau_1}, \omega^{(2)}_{\tau+})] \label{eqn:ST_future_cost}\\
\text{subject to:} \quad & \nonumber \\
& x_{g, \tau+} = x_{g, \tau} + v_{g, \tau} - w_{g, \tau} \label{eqn:ST_commit_dyn} \\
& p_{g, \tau+} = p_{g, \tau} + \Delta p_{g, \tau} \label{eqn:ST_power_dyn} \\
& -\Delta P_g^{\text{MAX}} \leq \Delta p_{g, \tau} \leq \Delta P_g^{\text{MAX}}  \label{eqn:ST_ramp_limit}\\
& (x_{g, \tau+}, p_{g, \tau+}) \in X^{(2)} \label{eqn:ST_state_feasible} \\
& y^{(2)}_{b, \tau} = y^{(1)}_{b} + \sum_{g \in \mathcal G^{(2)}} p_{g,\tau} \label{eqn:ST_agg_state_generation}
\end{align}

The term on line \eqref{eqn:ST_immediate_cost} models the immediate operational costs incurred at tick $\tau$. The compact notation $\mathbf{c}^{(2)}_\tau \mathbf{x}_\tau + \mathbf{d}^{(2)}_\tau \mathbf{u}_\tau$ expands to the sum of costs over all fast-ramping generators:
\[
\sum_{g \in \mathcal G^{(2)}} \left( C^f_g x_{g,\tau} + C^v_g p_{g,\tau} + C^{SU}_g v_{g,\tau} \right)
\]
where $C^f_g$ is the fixed (or no-load) cost for keeping a generator online, $C^v_g$ is the variable cost of power production, and $C^{SU}_g$ is the startup cost for generator $g$.

Line \eqref{eqn:ST_RT_cost} recursively includes the expected cost from the finer, timescale-3 (RT) subproblem that occurs between the current ST tick and the next. Since the RT timescale primarily handles the dispatch of renewable generators which lack complex state dynamics (e.g., commitment status), the value function $V^{(3)}_\tau$ does not require a prescribed ending state from the ST model. The final term, on line \eqref{eqn:ST_future_cost}, is the recursion along the ST timescale itself, representing the expected cost-to-go from the next tick, $\tau+$, onward, as is standard in a single timescale problem.

The conditioning in the expectation operators, reflects the hierarchical nature of the multi-timescale uncertainty. The first expectation, $E_{\omega^{(3)}_\tau | \omega^{(2)}_\tau}$, signifies that the probability distribution of the real-time (RT) randomness $\omega^{(3)}_\tau$ is dependent on the realized outcome of the short-term (ST) randomness $\omega^{(2)}_\tau$. The second expectation, $E_{\omega^{(2)}_{\tau+} | \omega^{(1)}_{\tau_1},\omega^{(2)}_\tau}$, describes the evolution of the ST process itself. The conditioning on $\omega^{(2)}_\tau$ reflects the Markovian property of the ST process, where the next random state depends on the current one. The conditioning on $\omega^{(1)}_{\tau_1}$, however, indicates that the transition probabilities of this entire Markovian process are parameterized by the day-ahead (DA) scenario. From the perspective of the ST subproblem, the DA randomness $\omega^{(1)}_{\tau_1}$ is treated as a given parameter that dictates the behavior of the underlying ST stochastic process for that segment.

The constraints \eqref{eqn:ST_commit_dyn} through \eqref{eqn:ST_state_feasible} define the system dynamics for the fast-ramping generators. Constraint \eqref{eqn:ST_commit_dyn} tracks the on/off commitment status based on startup and shutdown decisions. Constraints \eqref{eqn:ST_power_dyn} and \eqref{eqn:ST_ramp_limit} enforce the power output dynamics and the ramping limit. \eqref{eqn:ST_state_feasible} ensures that the resulting state of each generator remains within its feasible operating region, $X^{(2)}$.

After the dispatch for the fast-ramping generators is determined, Equation \eqref{eqn:ST_agg_state_generation} tallies the total power available at each bus $b$. It sums the newly dispatched power from the fast-ramping generators with the power already provided by the slower timescale-1 generators ($y^{(1)}_{b, \tau_1}$).

The dynamic programming recursion terminates at the end of each segment, specifically when the short-term tick $\tau$ go past the start of the next day-ahead tick, $\tau_1+$. The value function is defined by a terminal penalty that penalizes any deviation from the prescribed target state:
\begin{equation}
    Q^{(2)}_\tau(\mathbf{x}_\tau, \mathbf{\tilde x}, \dots) = \rho \| \mathbf{x}_{\tau} - \mathbf{\tilde x} \| \label{eqn:ST_terminal_cost}
\end{equation}

Finally, we define the value function $V^{(2)}_{\tau_1}$ that appears in the day-ahead objective function \eqref{eqn:DA_obj}. This function represents the total expected cost for the short-term segment corresponding to the day-ahead tick $\tau_1$. We define it in terms of our auxiliary function $Q^{(2)}$ evaluated at the very first tick $\tau_2$ of that segment:
\begin{equation}
    V^{(2)}_{\tau_1}(\mathbf{x}_{\tau_1}, \mathbf{\tilde x}, \mathbf{y}^{(1)}_{\tau_1};\omega^{(1)}_{\tau_1}) := Q^{(2)}_{\tau_1}(\mathbf{x}_{\tau_1}, \mathbf{\tilde x}, \mathbf{y}^{(1)}_{\tau_1};\omega^{(1)}_{\tau_1}, \omega^{(2),0}_{\tau_1}) \label{eqn:V_equals_Q}
\end{equation}
Here, $\omega^{(2),0}_{\tau_1}$ is a fixed, deterministic value used to initialize the ST Markovian random process.

\subsection{Real-Time Timescale}

The final and finest level of the hierarchy is the Real-Time (RT) timescale, whose expected operational cost is represented by the value function $V^{(3)}_\tau(y^{(2)}, \omega^{(2)}_\tau)$. This function solves an economic dispatch problem at each RT tick. Its goal is to meet the final load balance requirements by managing power transmission between buses, using the total available conventional power $y^{(2)}$ passed down from the ST model. This dispatch decision is made based on the realized outcome of the real-time renewable generation $\omega^{(3)}$, the distribution of which is conditioned on the ST forecast $\omega^{(2)}_\tau$. As this stage formulates a standard optimal power flow (OPF) problem, its detailed mathematical equations are well-documented in the literature and are omitted here for brevity.

\section{Conclusions} \label{section:conclusions}

We introduced a multi-timescale modeling framework that coordinates decisions between timescales using aggregated state variables and extended the ideas from the multi-horizon trees using synchronized states assumption. Given the rich structure of this problem, there can be a lot of interesting areas to explore. Here we list a few directions for future research:

\begin{itemize}
    \item We need to implement an algorithm to solve the resulting multi-timescale problem. The modeling tools from ``Plasmo.jl'' \cite{jalving_graph-based_2022} are very suitable for this application. In addition, we may be able to use problem-specific strategies to help scale up. We note that a lot of the subproblems in a multi-timescale problem are almost identical, so it may be possible to form heuristic solutions by collecting the solutions from another time segment in the same timescale or by forming an approximation of the faster timescales.
    \item Because of the sheer size of a practical multi-timescale problem, it is highly unlikely we can obtain a high quality solution through one single instantiation of the multi-timescale randomness. Incremental approaches, such as those described in \cite{zhang_sampling-based_2024}, that solves problem iteratively with more and more detailed description of the randomness, would be a promising improvement to tractability. Furthermore, scenario generation algorithms could selectively refine certain parts of the tree while maintaining coarser representations in less critical regions, guided by domain knowledge. For example, in power systems, wind dynamics during calm, non-peak periods may not require extensive modeling.
    \item On timescales that use value function instantiation, ensemble-based methods such as compromise solutions \cite{xu_compromise_2023} \cite{xu_ensemble_2024} could improve solution quality by using results from repetitive subproblem solves. This idea can be further enhanced with parametric/non-parametric representation of the randomness \cite{diao_distribution-free_2024} within our multi-timescale framework. Due to the similarity of the multistage subproblems in our multi-timescale problem, we expect these techniques will have a great impact.
\end{itemize}


\textbf{Acknowledgments} This research was supported by the Advanced Scientific Computing Research (ASCR) program under the U.S. Department of Energy, Award Number DE-SC0023361.

\begin{appendices}

\section{Multistage Stochastic Optimal Control} \label{section:appendix-multistage}

In a multistage stochastic optimal control problem, there is typically only one timescale. \footnote{Although the model could have multiple timescales, the timescales are not written out explicitly.} To define a multistage stochastic optimal control problem, we need the following elements:

\begin{itemize}
    \item Stages $\mathcal T = \{0, 1, 2, ..., T-1, T\}$
    \item State $x_t$
    \item Control $u_t \in U_t(x_t)$
    \item Randomness $\omega_t \in \Omega_t$
    \item State Dynamics $x_{t+1} = D_t(x_t, u_t, \omega_t)$
    \item Stage Cost $f_t(x_t, u_t)$ 
\end{itemize}

In a multistage problem, the set of stages $\mathcal T$ is typically comprised of integers. This set defines the sequence in which decisions will be made and serves as the index set for the variables. For each stage $t \in \mathcal T$, there is a state variable $x_t$, which describes the current status of the system we are interested in. We also have control variables $u_t$ that represent the actions that can be taken to influence the system. These controls must be feasible, meaning they must belong to the set $U_t(x_t)$, which defines the permissible actions based on the current state $x_t$.

For each stage \(t\), there is an associated random variable \(\omega_t \in \Omega_t\) that describes the uncertainty or unknown information, which is revealed to the decision maker after the control \(u_t\) is decided. Once the action \(u_t\) is taken and the uncertainty \(\omega_t\) is realized, the system transitions from state \(x_t\) to \(x_{t+1} = D_t(x_t, u_t, \omega_t)\). This transition incurs an instantaneous cost represented by \(f_t(x_t, u_t)\).

\textbf{It is important to understand \(x_t\) and \(u_t\) as decision random variables.} This is because \(x_{t+1}\) depends on \(\omega_t\) in the dynamics, making it a (decision) random variable as well. Furthermore, since \(u_t\) is selected from \(U_t(x_t)\) and \(x_t\) is a random variable, \(U_t(x_t)\) is actually a random set. Consequently, \(u_t\) must also be treated as a random variable.

\begin{figure}
    \centering
    \includegraphics[width=0.5\linewidth]{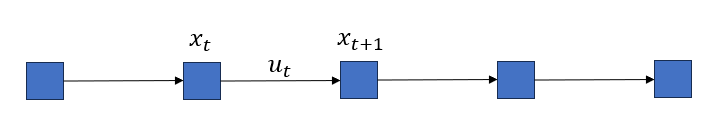}
    \caption{A visualization of a multistage problem (random variable view, each node is a state variable as a random decision variable)}
    \label{fig:multistage_problem}
\end{figure}

In this notation, we define a stochastic optimal control problem in mathematical terms. If we use blue squares (\img{square.PNG}) to denote each state, the representation can be visualized as shown in Figure \ref{fig:multistage_problem}.

\begin{align*}
    \text{Min} \ & E_{\omega} [\sum_{t \in \mathcal T} f_t(x_t, u_t)] \\
    & x_{t+1} = D_t(x_t, u_t, \omega_t) \\
    & u_t \in U_t(x_t) \\
    & x_t, u_t \text{ non-anticipative}
\end{align*}

\subsection{Non-anticipativity in single timescale problem} \label{section:nac_single_timescale}

One key requirement of stochastic programming is non-anticipativity, which means that decisions cannot be based on future unseen information. We can use a natural filtration $\mathcal F$ to describe the information available up to each point.

\[
\mathcal F: F_0 = \sigma(\{\}) \quad F_0 = \sigma(\omega_0)  \quad F_1 = \sigma(\omega_0, \omega_1) \quad  \dots
\]

In general, we define the information up to stage $t$ by the following sigma-algebra $F_t$. Specially, we define $F_0$ as the trivial sigma algebra because no information has been revealed so far. 

\begin{equation}
    F_t := \sigma(\{\tilde \omega_{t'}: t' \leq t\})
\end{equation}

To describe the non-anticipativity of $x_t$ and $u_t$, it means $x_t$ and $u_t$ are measurable with respect to $F_t$. Equivalently, we say they adapt to the filtration $\mathcal F$. It is sometimes notated as \(\{x_t, u_t\}_{t \in \mathcal T} \triangleleft \mathcal F\).

\begin{figure}
    \centering
    \includegraphics[width=\linewidth]{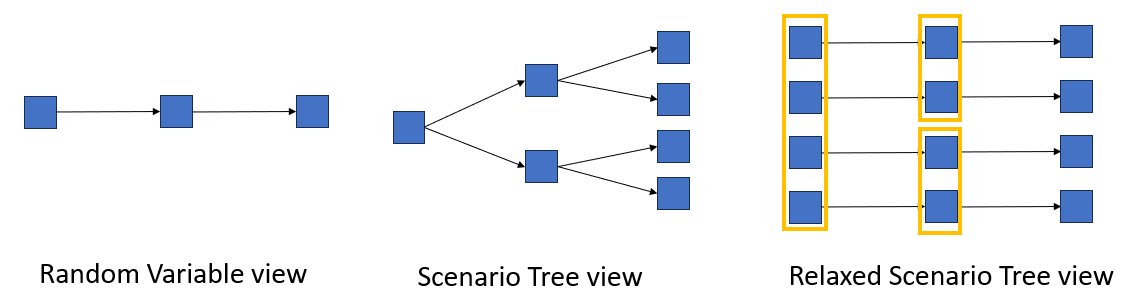}
    \caption{Different representations of a multistage problem}
    \label{fig:multistage_nac}
\end{figure}

We will use these visualizations to illustrate the problems, as shown in Figure \ref{fig:multistage_nac}, which describes the same multistage problem. In a random variable view, each node represents a decision random variable. In a scenario tree view, each node is understood as a decision vector in the usual sense. This representation can be derived by expanding the random variable view, assuming a set of scenarios to describe the random variables. In the relaxed scenario tree view, which is a further expansion of the scenario tree view, there are duplicate copies of the nodes. The nodes in each orange rectangle box should be equal to each other, representing non-anticipativity constraints.

\end{appendices}

\bibliography{references}

\end{document}